\documentclass[]{interact}

\usepackage{epstopdf}
\usepackage[caption=false]{subfig}

\usepackage[numbers,sort&compress]{natbib}
\bibpunct[, ]{[}{]}{,}{n}{,}{,}
\makeatletter
\def\NAT@def@citea{\def\@citea{\NAT@separator}}
\makeatother
\usepackage{anyfontsize}
\theoremstyle{plain}
\newtheorem{theorem}{Theorem}[section]
\newtheorem{lemma}[theorem]{Lemma}
\newtheorem{corollary}[theorem]{Corollary}
\newtheorem{proposition}[theorem]{Proposition}

\theoremstyle{definition}
\newtheorem{definition}[theorem]{Definition}
\newtheorem{remark}[theorem]{Remark}

\usepackage{xcolor}

\theoremstyle{remark}

\begin{document}
	

	\title{Hyers-Ulam stability of closed linear relations in Hilbert spaces}
 \author{\name {Arup Majumdar \thanks{Arup Majumdar (corresponding author). Email address: arupmajumdar93@gmail.com}}\affil{Department of Mathematical and Computational Sciences, \\
			National Institute of Technology Karnataka, Surathkal, Mangaluru 575025, India.}}
  
	\maketitle

	\begin{abstract}
		This paper introduces the concept of Hyers-Ulam stability for linear relations in normed linear spaces and presents several intriguing results that characterize the Hyers-Ulam stability of closed linear relations in Hilbert spaces. Additionally, sufficient conditions are established under which the sum and product of two Hyers-Ulam stable linear relations remain stable. 
    \end{abstract}
	
      \begin{amscode}47A06; 47A10; 47B25.\end{amscode}
    \section{Introduction}
   The theory of Hyers-Ulam stability plays a significant role in various branches of mathematics, including functional equations, optimization, differential equations, and statistics. The stability problem for functional equations was first introduced by Ulam in a 1940 talk at the University of Wisconsin, where he posed the question: ``For what metric groups $G$, is it true that an $\varepsilon$-automorphism of $G$ is necessarily close to an automorphism?" The question was answered by Hyers in 1941 for Banach spaces: Let $X$ and $Y$ be two real Banach spaces and $f: X \to Y$ be a mapping such that for each fixed $x \in X$, $f(tx)$ is continuous in $t\in \mathbb{R}$ (the set of all real numbers), and if there exists  $\varepsilon \geq 0$ satisfying the inequality $$\| f(x+y) - f(x) -f(y)\| \leq \varepsilon \quad \text{ for all }x, y \in X,$$ then there exists a unique linear mapping $L: X \to Y$ such that $\|f(x) - L(x)\| \leq \varepsilon$ for every $x \in X$. This result is known as the Hyers-Ulam stability of the additive Cauchy equation $$g(x+y) = g(x) + g(y).$$
  In 1978, Rassias generalized this by considering an unbounded right-hand side in the involved inequalities, depending on certain functions, which led to the concept of modified Hyers-Ulam stability for the additive functional equation \cite{MR1778016}. Following this, numerous researchers extended Ulam's stability problem to other types of functional equations, leading to generalizations of Hyers' results in various directions.
   
   Obloza \cite{MR1321558} was the first to establish results on the Hyers-Ulam stability of differential equations. Alsina and Ger \cite{MR1671909} explored the Hyers-Ulam stability for first-order linear differential equations. Miura et al. further generalized the results for $n^\text{th}$ order linear differential operator $p(D)$ and proved that the differential operator equation $$p(D)f=0$$ is Hyers-Ulam stable if and only if the algebraic equation $p(z)= 0$ has no pure imaginary solution, where $p$ is a complex-valued polynomial of degree $n$, and $D$ is a differential operator \cite{MR2000046}.  In the same paper, Miura et al. first introduced the concept of the Hyers-Ulam stability of a mapping (not necessarily linear) between two complex linear spaces $X$ and $Y$ with gauge functions $\rho_X$ and $\rho_Y$, respectively. A mapping $S$ has the Hyers-Ulam stability (HUS) if there exists a constant $M \geq 0$ with the following property:
	For every $\varepsilon \geq 0$, $y\in S(X)$ and $x\in X$ satisfying $\rho_Y (S(x) -y) \leq \varepsilon$ we can find an $x_0 \in X$ such that $S(x_0) =y$ and $\rho_X (x-x_0) \leq M \varepsilon,$ where $M$ is called as a HUS constant, and the infimum of all the HUS constants for $S$ by $M_{S}$. Essentially, if $S$ has the HUS, then for each $y\in S(X)$ and ``$\varepsilon$-approximate solution" $x$ of the equation $S(u)= y$ there corresponds an exact solution $x_{0}$ of the equation that is contained in a $M\varepsilon$- neighbourhood of $x$. Subsequently, Hirasawa and Miura expanded the concept of Hyers-Ulam stability of closed operators in Hilbert spaces in 2006 \cite{MR2204863}. Moreover, In the same paper \cite{MR2204863}, they established the Hyers-Ulam stability of a linear operator $T$ from the domain $D(T) \subset X$ into $Y$, where $X$ and $Y$ both are normed linear spaces. Specifically, there exists a constant $M \geq 0$ with the following property:\\
		For any $y\in R(T), ~ \varepsilon \geq 0$ and $x\in D(T)$ with $\| Tx -y\| \leq \varepsilon$, there exists $x_0 \in D(T)$ such  that $Tx_{0} = y$ and $\|x -x_{0}\| \leq M\varepsilon$.\\
    In 2024, Majumdar et al. investigated the Hyers-Ulam stability of closable operators in Hilbert spaces and provided several characterizations of Hyers-Ulam stable closable operators \cite{majumdar2024hyers}. This paper delves into the exploration of the Hyer-Ulam stability of closed linear relations in Hilbert spaces. Section 2 is dedicated to the basic definitions and notations related to linear relations. In Section 3, we discuss several properties of the Hyers-Ulam stable linear relations in Hilbert spaces.
 \section{Preliminaries}
 Throughout the paper, the symbols $H, K, H_{i}, K_{i} ~(i = 1, 2)$ represent real or complex Hilbert spaces. A linear relation $T$ from $H$ into $K$ is a linear subspace of the Cartesian product in $H \times K$, and the collection of all linear relations from $H$ into $K$ is denoted by $LR(H, K)$. We call $T$ a closed linear relation from $H$ into $K$ if it is a closed subspace of $H \times K$ and the set of all closed linear relations from $H$ into $K$ is denoted by $CR(H, K)$. The following notations of domain, range, kernel and multi-valued part of a linear relation $T$ from $H$ into $K$ will be used respectively in the paper:\\
 $$D(T) = \{h \in H: \{h, k\} \in T\},  ~ R(T) = \{k \in K: \{h, k\} \in T\}$$
 $$ N(T)= \{h\in H:\{h, 0\} \in T\}, ~ M(T) = \{k\in K: \{0, k\} \in T \}.$$
 It is obvious that $N(T)$ and $M(T)$ both are closed subspaces in $H$ and $K$, respectively, whenever $T$ is a closed linear relation from $H$ into $K$. In general, the inverse of an operator is always a linear relation. The inverse of a linear relation $T$ from $H$ into $K$ is defined as $T^{-1}= \{\{k, h\} \in K \times H: \{h, k\} \in T \}$. Thus, it is immediate that $D(T^{-1}) = R(T)$ and $N(T^{-1}) = M(T)$. We define $Tx = \{y \in K: \{x, y\} \in T\}, \text{ where } T \in L(H, K).$ Consequently, $T\vert_{W}$ denotes the restriction of $T \in LR(H, K)$ with domain $D(T) \cap W$, where $W$ is a subset of $H$ (in other words, $T\vert_{W}$ is equal to $T$ in domain $D(T) \cap W$). The adjoint of a linear relation $T$ from $H$ into $K$ is the closed linear relation $T^{*}$ from $K$ into $H$ defined by:
 $$T^{*} = \{\{k^{'}, h^{'}\} \in K \times H : \langle h^{'}, h \rangle = \langle k^{'}, k \rangle, \text{ for all } \{h, k\} \in T\}.$$
 Observe that $(T^{-1})^{*} = (T^{*})^{-1}$, so that $(D(T))^{\perp} = M(T^{*})$ and $N(T^{*}) = (R(T))^{\perp}$. A linear relation $T$ in $H$ is said to be symmetric if $T \subset T^{*}$. Again, a linear relation $T$ in $H$ is non-negative if $\langle k, h \rangle \geq 0$, for all $\{h, k\} \in T$. Moreover, a linear relation $T$ in $H$ is said to be self-adjoint when $T = T^{*}$. If $S$ and $T$ both are linear relations then their product $TS$ defined by:
 $$TS = \{\{x,y\}: \{x,z\} \in S \text{ and } \{z,y\} \in T, \text{ for some } z \}.$$
 The sum of two linear relations $T$ and $S$ from $H$ into $K$ is $T + S = \{\{x, y +z\}: \{x, y\} \in T \text{ and } \{x, z\} \in S \}$. Whereas the Minkowski sum is denoted by $T \widehat{+} S:=\{\{x+v, y+w\}:\{x, y\} \in T,\{v, w\} \in S\} .$ \\
 Here, $Q_{T}$ denotes the natural quotient map from $K$ into $K / \overline{M(T)}$, where $T$ is a linear relation from $H$ into $K$. It is easy to show that $Q_{T}T$ is a linear operator from $H$ into $K / \overline{M(T)}$. We call the linear relation $T$ from $H$ into $K$ bounded linear relation if $Q_{T}T$ is a bounded operator and the set of all bounded linear relations from $H$ into $K$ is denoted by $BR(H, K)$. Some characterizations of bounded linear relations are explored in \cite{MR1631548}.\\
 The regular part of a closed linear relation $T$ from $H$ into $K$ is $P_{\overline{D(T^{*})}}T$, denoted by $T_{op}$ which is an operator with $T_{op} \subset T$, where $P_{\overline{D(T^{*})}}$ is the orthogonal projection in $K$ onto $\overline{D(T^{*})} = (M(T))^{\perp}$. It can be shown that $T = T_{op} \widehat{+} ~(\{0\} \times M(T))$, when $T$ is a closed linear relation from $H$ into $K$ \cite{MR3971207}. When $T$ is a closed relation from $H$ into $K$, then $T_{op}$ is also a closed operator \cite{MR3971207}. Finally, the linear operator $(T^{-1})_{op} = P_{(N(T))^{\perp}}T^{-1}$ is called the Moore-Penrose inverse of the linear relation $T$ from $H$ into $K$, denoted by $T^{\dagger}$.
\section{Characterizations of the Hyers-Ulam stable closed linear relations in Hilbert spaces} 
\begin{definition}\label{def 3.1}
Let $T$ be a linear relation from a normed linear space $X$ into a normed linear space $Y$. Then $T$ is said to be the Hyers-Ulam stable if there exists a constant $M \geq 0$ with the following property:
for any $y \in R(T)$, $\varepsilon \geq 0$, and $y_{0}\in R(T)$ with $\|y-y_{0}\| \leq \varepsilon$, there exist $\{x, y\} \in T$ and $\{x_{0}, y_{0}\} \in T$ such that $\|x-x_{0}\| \leq M \varepsilon$.
\end{definition}
We call $M$ a Hyers-Ulam stable (HUS) constant for the linear relation $T$, and the infimum of all HUS constants of $T$ is denoted by $M_{T}$.
\begin{remark}\label{remark 3.2}
Let $T$ be a linear relation from a normed linear space $X$ into a normed linear space $Y$. If $T$ is Hyers-Ulam stable, then there exists a constant $M \geq 0$ with the following property:
for any $y \in R(T)$ and $y_{0} \in R(T)$, there exist $\{x, y\} \in T$ and $\{x_{0}, y_{0}\} \in T$ such that $\|x-x_{0}\| \leq M \|y -y_{0}\|$.
\end{remark}
From now on, we consider $T$ to be a closed linear relation from $H$ into $K$.
\begin{theorem}\label{thm 3.3}
Let $T \in CR(H, K)$. Then $T$ is Hyers-Ulam stable if and only if $R(T)$ is closed.
\end{theorem}
\begin{proof}
Since $T$ is closed implies $T^{-1}$ is also closed. I claim that $T^{-1}$ is bounded in order to show that $D(T^{-1}) = R(T)$ is closed (by Theorem III.4.2 \cite{MR1631548}). Let us consider $y \in R(T)$ and $y_{0} \in R(T)$, there exist $\{x, y\} \in T$ and $\{x_{0}, y_{0}\} \in T$ such that $\|x-x_{0}\| \leq M\|y -y_{0}\|$, where $M$ is a HUS constant of $T$. Then, $\{x-x_{0}, y -y_{0}\} \in T$. Thus, $\|Q_{T^{-1}}T^{-1}(y -y_{0})\| = \|(x -x_{0}) + N(T)\| \leq \|x -x_{0}\| \leq M \|y-y_{0}\|$. So, $Q_{T^{-1}}T^{-1}$ is bounded implies $T^{-1}$ is bounded. Hence, $R(T) = D(T^{-1})$ is closed.\\
Conversely, suppose $R(T)$ is closed. Then $Q_{T^{-1}}T^{-1}$ is bounded because $T^{-1}$ is bounded. Then for $z \in R(T)$ and $z_{0} \in R(T)$, we have $\{u, z\} \in T$ and $\{u_{0}, z_{0}\} \in T$, for some $u, u_{0} \in D(T)$. Now,
$$\|(u-u_{0}) + N(T)\| = \|Q_{T^{-1}}T^{-1} (z -z_{0})\| \leq \|T^{-1}\| \|z-z_{0}\|.$$
We get $v \in N(T)$ such that $\|u-u_{0}+ v\| \leq (\|T^{-1}\|+1)\|z-z_{0}\|$. Moreover, $\{u, z\} \in T$ and $\{u_{0}-v, z_{0}\} \in T$. Therefore, $T$ is Hyers-Ulam stable. It is obvious to show that $M_{T} = \|T^{-1}\|.$
\end{proof}
\begin{theorem}\label{thm 3.4}
    Let $T \in CR(H, K)$ be Hyers-Ulam stable. Then $M_{T}$ is a Hyers-Ulam stable constant. 
\end{theorem}
\begin{proof}
First, we assume that $M_{T}$ is not a Hyers-Ulam stable constant. Then there exist $y \in R(T)$ and $y_{0} \in R(T)$ such that for all $x \in T^{-1}y$ and $x_{0} \in T^{-1}y_{0}$, we have $\|x-x_{0}\| > M_{T}\|y-y_{0}\|.$\\
It is easy to show that $T^{-1}y$, and $T^{-1}y_{0}$ both are non-empty closed convex subsets in $H$. By Corollary 16.6(a) \cite{MR1427262} (which states that every nonempty convex subset of Hilbert spaces contains an element of minimal norm), we get
\begin{align*}
dist(T^{-1}y, T^{-1}y_{0}) &= \inf_{z_{0} \in T^{-1}y_{0}} \inf_{z \in T^{-1}y}\|z_{0}-z\| \\
&= \inf_{z_{0} \in T^{-1}y_{0}} \|z_{0}-z^{'}\|, \text{ for some }z^{'} \in T^{-1}y \\
&= \|z_{0}^{'} - z^{'}\|, \text{ for some } z_{0}^{'} \in T^{-1}y_{0}.
\end{align*}
Thus, $dist(T^{-1}y, T^{-1}y_{0}) = \|z^{'}-z_{0}^{'}\| > M_{T}\|y-y_{0}\|$. There exists a positive $\delta > 0$ such that $M_{T} + \delta$ is a Hyers-Ulam stable constant and $dist(T^{-1}y, T^{-1}y_{0})  > (M_{T} + \delta)\|y-y_{0}\|$. Hence, for all $x \in T^{-1}y$ and $x_{0} \in T^{-1}y_{0}$, we have $\|x-x_{0}\| > (M_{T}+ \delta) \|y-y_{0}\|$, which is a contradiction. Therefore, $M_{T}$ is a Hyers-Ulam stable constant.
\end{proof}
\begin{theorem}\label{thm 3.5}
Let $T \in CR(H, K)$. Then $T$ is Hyers-Ulam stable if and only if $T_{op}$ is Hyers-Ulam stable.
\end{theorem}
\begin{proof}
Suppose $T\in CR(H, K)$ is Hyers-Ulam stable. By Theorem \ref{thm 3.3}, $R(T)$ is closed. We claim $R(T_{op})$ is closed in order to show the Hyers-Ulam stability of $T_{op}$. Let $y \in \overline{R(T_{op})}$, then there exists $\{y_{n}\}$ in $R(T_{op})$ such that $y_{n} \to y$, as $n \to \infty$. We know, $T_{op} = P_{\overline{D(T^{*})}}T$, so there exist $\{x_{n}, z_{n}\} \in T$ and $\{z_{n}, y_{n}\} = \{z_{n}, P_{\overline{D(T^{*})}}z_{n}\} \in G(P_{\overline{D(T^{*})}})$ (for all $n \in \mathbb{N}$), where $G(P_{\overline{D(T^{*})}})$ is the graph of operator $P_{\overline{D(T^{*})}}$. We can write $z_{n} = y_{n} + y_{n}^{'}$, where $y_{n}^{'} \in M(T)$ for all $n \in \mathbb{N}$. It confirms that $\{x_{n}, y_{n}\} = \{x_{n}, z_{n}-y_{n}^{'}\} \in T$. Thus, $y \in \overline{D(T^{*})} \cap R(T).$ There exists a $x \in D(T)$ such that $\{x, y\} \in T$ and $\{y, y\} \in G(P_{\overline{D(T^{*})}})$ which implies $\{x, y\} \in T_{op}$. Hence, $R(T_{op})$ is closed and $T_{op}$ is Hyers-Ulam stable.\\
Conversely, $T_{op}$ is Hyers-Ulam stable. we will show that $R(T)$ is closed to prove the Hyers-Ulam stability of $T$. Consider $w \in \overline{R(T)}$, there exists a sequence $\{w_{n}\}$ in $R(T)$ such that $w_{n} \to w$, as $n \to \infty$. So, there exists a sequence $\{v_{n}\}$ in $D(T)$ such that $\{v_{n}, w_{n}\} \in T = T_{op} \widehat{+} (\{0\} \times M(T))$, for all $n \in \mathbb{N}$. Thus, for all $n \in \mathbb{N}$, $\{v_{n}, w_{n}\} = \{v_{n}, w_{n}^{'}\} + \{0, w_{n}^{''}\}$, where $\{v_{n}, w_{n}^{'}\} \in T_{op}$ and $\{0, w_{n}^{''}\} \in (\{0\} \times M(T))$. Again, the convergent sequence $\{w_{n}\}$ says that $\{w_{n}^{'}\}$ and $\{w_{n}^{''}\}$ both are convergent to $w^{'}$ and $w^{''}$ respectively, for some $w^{'} \in R(T_{op}) \subset (M(T))^{\perp}$ and $w^{''} \in M(T)$. This guarantees that $w = w^{'} + w^{''}$. There is $u \in D(T_{op})$ such that $\{u, w^{'}\} \in T_{op}$ and $\{0, w^{''}\} \in (\{0\} \times M(T))$ which implies $\{u, w\} \in T$ and $w \in R(T)$. Therefore, $T$ is Hyers-Ulam stable.
\end{proof}
\begin{proposition}\label{pro 3.6}
Let $T \in CR(H, K)$. Then $T$ is Hyers-Ulam stable if and only if $T^{*}$ is Hyers-Ulam stable.
\end{proposition}
\begin{proof}
Since, $T$ is in $CR(H, K)$. By Proposition 2.5 \cite{MR3079830} and Theorem \ref{thm 3.3}, we get $T$ is Hyers-Ulam stable if and only if $R(T)$ is closed if and only if $R(T^{*})$ is closed if and only if $T^{*}$ is Hyers-Ulam stable.
\end{proof}
\begin{proposition}\label{pro 3.7}
Let $T \in CR(H, K)$. Then $T$ is Hyers-Ulam stable if and only if $T^{*}T$ is Hyers-Ulam stable ($TT^{*}$ is Hyers-Ulam stable).
\end{proposition}
\begin{proof}
By Lemma 5.1 \cite{MR2188974}, we get $T^{*}T$ and $TT^{*}$ both are non-negative self-adjoint. So, $T^{*}T$ and $TT^{*}$ both are closed. By Proposition 2.5 \cite{MR3079830} and Theorem \ref{thm 3.3}, we have $T$ is Hyers-Ulam stable if and only if $R(T)$ is closed if and only if $R(T^{*})$ is closed if and only if $R(T^{*}T)$ is closed ($R(TT^{*})$ is closed) if and only if $T^{*}T$ is Hyers-Ulam stable ($TT^{*}$ is Hyers-Ulam stable).
\end{proof}
\begin{theorem}\label{thm 3.8}
Let $T \in CR(H, K)$. $T^{\dagger}$ is Hyers-Ulam stable if and only if $T$ is bounded.
\end{theorem}
\begin{proof}
Theorem III.4.2 \cite{MR1631548} says that $T$ is bounded if and only if $D(T)$ is closed. Since, $T$ is closed, so $T^{-1}$ is closed implies $T^{\dagger} = (T^{-1})_{op}$ is closed. By Theorem \ref{thm 3.5} and Theorem \ref{thm 3.3}, we get $D(T)= R(T^{-1})$ is closed. Hence, $T$ is bounded.\\ 
Conversely, suppose $T$ is bounded. So, $D(T) = R(T^{-1})$ is closed. Then, $T^{-1}$ is Hyers-Ulam stable. By Theorem \ref{thm 3.5}, we have the Hyers-Ulam stability of $T^{\dagger} = (T^{-1})_{op}$. 
\end{proof}
Let $T \in CR(H, K)$ be a non-negative self-adjoint linear relation. Then $T_{op}$ is also a non-negative self-adjoint linear operator \cite{MR3971207}. Moreover, $T^{\frac{1}{2}} = (T_{op})^{\frac{1}{2}} \widehat{+} (\{0\} \times M(T))$ and $(T_{op})^{\frac{1}{2}} = (T^{\frac{1}{2}})_{op}$ \cite{MR3971207}.
\begin{lemma}\label{lemma 3.9}
Let $T \in CR(H, K)$ be a non-negative self-adjoint linear relation. $T$ is Hyers-Ulam stable if and only if $T^{\frac{1}{2}}$ is Hyers-Ulam stable. 
\end{lemma}
\begin{proof}
The self-adjointness of $T_{op}$ says that $T_{op}$ is densely defined in $\overline{D(T)} = \overline{D(T^{*})}$ (by Theorem 1.3.16 \cite{MR3971207}). By Theorem \ref{thm 3.5} we get, $T$ is Hyers-Ulam stable if and only if $T_{op}$ is Hyers-Ulam stable if and only if $ (T^{\frac{1}{2}})_{op} = (T_{op})^{\frac{1}{2}}$ is Hyers-Ulam stable (by Proposition 2.23 \cite{majumdar2024hyers}) if and only if $T^{\frac{1}{2}}$ is Hyers-Ulam stable (by Theorem \ref{thm 3.5}).
\end{proof}
\begin{theorem}\label{thm 3.10}
Let $T \in CR(H, K)$. $T$ is Hyers-Ulam stable if and only if $\vert T \vert = (T^{*}T)^{\frac{1}{2}}$ is Hyers-Ulam stable. 
\end{theorem}
\begin{proof}
By Proposition \ref{pro 3.7} and Lemma \ref{lemma 3.9}, we get $T$ is Hyers-Ulam stable if and only if $T^{*}T$ is Hyers-Ulam stable if and only if $\vert T \vert$ is Hyers-Ulam stable.
\end{proof}
\begin{theorem}\label{thm 3.11}
Let $T \in CR(H, K)$. Then $C_{T} = (I + T^{*}T)^{-1}$ is Hyers-Ulam stable if and only if $T$ is bounded.
\end{theorem}
\begin{proof}
Theorem 5.2 \cite{MR4203636} says that $C_{T} = P_{T}P_{T}^{*}$, where $P_{T}\{x, y\} = x$, for all $\{x, y\} \in T$. So, $C_{T}$ is a bounded operator in domain $H$, which implies $C_{T}$ is closed. First consider $C_{T}$ is Hyers-Ulam stable then $R(C_{T})= R((I + T^{*}T)^{-1}) = D(T^{*}T)$ is closed. By Lemma 5.1 (a) \cite{MR4203636} we get, $D(T^{*}T) = \overline{D(T)} \subset D(T) \subset \overline{D(T)}$ implies $D(T)$ is closed. Hence, $T$ is bounded (by Theorem III.4.2 \cite{MR1631548}).\\
Conversely, suppose $T$ is bounded. Then, $T_{op}$ is bounded because $T$ is bounded and $\|T_{op}x\|  = \|Tx + M(T)\| = \|Tx\|$, for all $x \in D(T_{op})$. By Corollary III.1.13 \cite{MR1631548}, $(T_{op})^{*}$ is a bounded linear relation. Then $\|(T_{op})^{*}T_{op}x\| \leq \|(T_{op})^{*}\|\|T_{op}\| \|x\|$. Thus, $(T_{op})^{*}T_{op}$ is bounded. By Lemma 5.1(b) \cite{MR4203636}, we get $T^{*}T$ is bounded which implies $D(T^{*}T)$ is closed. Thus, $R(C_{T})$ is closed. Therefore, $C_{T}$ is Hyers-Ulam stable.
\end{proof}
\begin{theorem}\label{thm 3.12}
Let $T \in BR(H, K) \cap CR(H, K)$. If $T$ is Hyers-Ulam stable, then $Z_{T} = T(I + T^{*}T)^{-\frac{1}{2}}$ is also Hyers-Ulam stable.
\end{theorem}
\begin{proof}
Since $C_{T}$ is a non-negative self-adjoint bounded operator in domain $H$. So, $C_{T}^{\frac{1}{2}}$ exists and $Z_{T}= TC_{T}^{\frac{1}{2}}$. Now, we claim that $Z_{T}$ is closed. Consider $\{x, y\} \in \overline{Z_{T}}$, then there exists a sequence $\{\{x_{n}, y_{n}\}\}$ in $Z_{T}$ such that $\{x_{n}, y_{n}\} \to \{x,y\}$ as $n \to \infty$, where $\{x_{n}, C_{T}^{\frac{1}{2}}x_{n}\} \in G(C_{T}^{\frac{1}{2}})$, and $\{C_{T}^{\frac{1}{2}}x_{n}, y_{n}\} \in T$, for all $n \in \mathbb{N}$. Since, $T$ is closed and $C_{T}^{\frac{1}{2}}$ is bounded in domain $H$. Thus, $\{x, y\} \in Z_{T}$ and $Z_{T}$ is closed. Moreover,
\begin{align*}
C_{T}^{\frac{1}{2}} = (I + T^{*}T)^{-\frac{1}{2}} = (I + (T_{op})^{*}T_{op})^{-\frac{1}{2}} = C_{T_{op}}^{\frac{1}{2}}.
\end{align*}
Again, $M(TC_{T}^{\frac{1}{2}}) = M(T)$ says that $\overline{D((TC_{T}^{\frac{1}{2}})^{*})} = \overline{D(T^{*})}$. Thus, 
\begin{align}\label{equ 1}
(Z_{T})_{op} = P_{\overline{D((TC_{T}^{\frac{1}{2}})^{*})}} (TC_{T}^{\frac{1}{2}}) = (P_{\overline{D(T^{*})}} T) C_{T}^{\frac{1}{2}} = T_{op}C_{T_{op}}^{\frac{1}{2}} = Z_{T_{op}}.
\end{align}
By Theorem 1.3.15 \cite{MR3971207}, we see that $Z_{T_{op}}$ is closed. From Theorem \ref{thm 3.5}, it suffices to show that $Z_{T_{op}}$ is Hyers-Ulam stable or that $R(Z_{T_{op}})$ is closed. Let $0 \neq w \in \overline{R(Z_{T_{op}})}$, then there exists a sequence $\{u_{n}\}$ in $H$ such that $T_{op}(P_{T}P_{T}^{*})^{\frac{1}{2}}u_{n} \to w$, as $n \to \infty$. Again, $R(T_{op})$ is closed because $T_{op}$ is Hyers-Ulam stable (by Theorem \ref{thm 3.5}). So, there exists an element $z \in D(T) \cap N(T_{op})^{\perp}$ such that $T_{op}z = w$ and
\begin{align}\label{equ 2}
\|T_{op}((P_{T}P_{T}^{*})^{\frac{1}{2}} u_{n}-z)\| \to 0, \text{ as } n \to  \infty.
\end{align}
 Let $v \in N(T)$, then $\{v, 0\} = \{v_{1}, v_{2}\} + \{0, v_{3}\}$, where $\{v_{1}, v_{2}\} \in T_{op}$ and $\{0, v_{3}\} \in (\{0\} \times M(T))$. This confirms $v_{1} = v$ and $v_{2} = v_{3} = 0$. Thus, $N(T) \subset N(T_{op}) \subset N(T)$ which implies $N(T) = N(T_{op})$. Moreover, $\gamma(T_{op}) > 0$ ($\gamma(T_{op})$ is the reduced minimum modulus of $T_{op}$) because $R(T_{op})$ is closed. By the relation (\ref{equ 2}), we get 
 \begin{align*}
 \gamma(T_{op}) \|P_{(N(T))^{\perp}}((P_{T}P_{T}^{*})^{\frac{1}{2}} u_{n}) -z\| \leq \|T_{op}((P_{T}P_{T}^{*})^{\frac{1}{2}} u_{n}-z)\| \to 0, \text{ as } n \to  \infty.
 \end{align*}
 Thus, $P_{(N(T))^{\perp}}((P_{T}P_{T}^{*})^{\frac{1}{2}} u_{n}) \to z$ as $n \to \infty$. Again, $R((P_{T}P_{T}^{*})^{\frac{1}{2}}) = R(P_{T}P_{T}^{*}) = R(P_{T}) = D(T)$ because $R(P_{T}) = D(T)$ is closed, since $T \in BR(H, K) \cap CR(H, K)$. Furthermore, $R((P_{T}P_{T}^{*})^{\frac{1}{2}}) + N(P_{(N(T))^{\perp}}) = D(T) + N(T) = D(T)$ is closed.  By Corollary 6 \cite{MR888141}, we have $R(P_{(N(T))^{\perp}}(P_{T}P_{T}^{*})^{\frac{1}{2}})$ is closed. Hence, there exists $q \in H$ such that $z = P_{(N(T))^{\perp}}((P_{T}P_{T}^{*})^{\frac{1}{2}})q$ and $w = T_{op}P_{(N(T))^{\perp}}((P_{T}P_{T}^{*})^{\frac{1}{2}})q = T_{op}((P_{T}P_{T}^{*})^{\frac{1}{2}})q$ (because, $N(T_{op}) = N(T)$). Now, we can say that $R(Z_{T_{op}})$ is closed which implies $(Z_{T})_{op}$ is Hyers-Ulam stable. Therefore, $Z_{T}$ is Hyers-Ulam stable.
\end{proof}
\begin{remark}\label{remark 3.13}
Let $T \in CR(H, K)$. When $T^{-1}$ is Hyers-Ulam stable if and only if $T$ is bounded. Because $T^{-1}$ is closed, and $T^{-1}$ is Hyers-Ulam stable if and only if $D(T) = R(T^{-1})$ is closed (by Theorem \ref{thm 3.3}) if and only if $T$ is bounded (by Theorem III.4.2 \cite{MR1631548}).
\end{remark}
\begin{lemma}\label{lemma 3.14}
Let $T \in CR(H, K)$. Then $T^{*}T = T^{*}T\vert_{(N(T))^{\perp}} \widehat{+} ~T^{*}T\vert_{N(T)}.$ Moreover, $T^{*}T\vert_{(N(T))^{\perp}}$ and $T^{*}T\vert_{N(T)}$ both are closed.
\end{lemma}
\begin{proof}
It is obvious to show that $T^{*}T \supset T^{*}T\vert_{(N(T))^{\perp}} \widehat{+} ~T^{*}T\vert_{N(T)}.$ Now consider $\{x, y\} \in T^{*}T$, then $x = x_{1} + x_{2} \in D(T) = N(T) + (N(T))^{\perp} \cap D(T)$, where $x_{1} \in N(T)$ and $x_{2} \in (N(T))^{\perp} \cap D(T)$. Then $\{x_{1}, 0\} \in T^{*}T\vert_{N(T)}$. There exists $z \in K$ such that $\{x, z\} \in T$ and $\{z, y\} \in T^{*}$. So, $\{x_{2}, z\} = \{x, z\} - \{x_{1}, 0\} \in T$. Thus, $\{x_{2}, y\} \in T^{*}T\vert_{(N(T))^{\perp}}$. This guarantees the reverse inclusion $T^{*}T \subset T^{*}T\vert_{(N(T))^{\perp}} \widehat{+} ~T^{*}T\vert_{N(T)}.$ Hence, $T^{*}T = T^{*}T\vert_{(N(T))^{\perp}} \widehat{+} ~T^{*}T\vert_{N(T)}.$\\
Now, we claim that $T^{*}T\vert_{N(T)}$ and  $T^{*}T\vert_{(N(T))^{\perp}}$ both are closed. Let $\{u, v\} \in \overline{T^{*}T\vert_{N(T)}}$. Then there exists a sequence $\{\{u_{n}, v_{n}\}\}$ in $T^{*}T\vert_{N(T)}$ with $u_{n} \in N(T)$ (for all $n \in \mathbb{N}$) such that $\{u_{n}, v_{n}\} \to \{u,v\}$ as $n \to \infty$. Again, $\{u_{n}, 0\} \in T$ and $\{0, v_{n}\} \in T^{*}$, for all $n \in \mathbb{N}$. The closeness of $T$, $T^{*}$ and $N(T)$ confirm that $\{u, v\} \in T^{*}T\vert_{N(T)}$. Thus, $T^{*}T\vert_{N(T)}$ is closed.\\
Consider, $\{s,t\} \in \overline{T^{*}T\vert_{(N(T))^{\perp}}}$, then there exists a sequence $\{\{s_{n}, t_{n}\}\}$ in $T^{*}T\vert_{(N(T))^{\perp}}$ with $s_{n} \in (N(T))^{\perp}$ (for all $n \in \mathbb{N}$) such that $\{s_{n}, t_{n}\} \to \{s,t\}$ as $n \to \infty$. So, $s \in (N(T))^{\perp}$. Again, $\{s,t\} \in T^{*}T$ because $T^{*}T$ is closed. Thus, $s \in D(T) \cap (N(T))^{\perp}$ and $\{s, t\} \in T^{*}T\vert_{(N(T))^{\perp}}$. Therefore, $T^{*}T\vert_{(N(T))^{\perp}}$ is closed.
\end{proof}
\begin{theorem}\label{thm 3.15}
Let $T \in CR(H, K)$. Then $\sigma(T^{*}T\vert_{(N(T))^{\perp}}) \setminus \{0\} = \sigma(T^{*}T) \setminus \{0\}.$
\end{theorem}
\begin{proof}
First, we claim $\lambda \in \sigma(T^{*}T) \setminus \{0\}$ implies $\lambda \in \sigma(T^{*}T\vert_{(N(T))^{\perp}}) \setminus \{0\}.$ Assume that $\lambda \in \rho(T^{*}T\vert_{(N(T))^{\perp}})$, where $T^{*}T\vert_{(N(T))^{\perp}}$ is a linear relation from Hilbert space $(N(T))^{\perp}$ into $(N(T))^{\perp}$. Then $(T^{*}T\vert_{(N(T))^{\perp}} - \lambda)^{-1} = ((T^{*}T -\lambda)\vert_{(N(T))^{\perp}}))^{-1}$ is bounded operator in domain $(N(T))^{\perp}$. Consider, $\{0, p\} \in (T^{*}T - \lambda)^{-1}$ implies $\{p, 0\} \in (T^{*}T - \lambda)$. Then $\{p, \lambda p\} \in T^{*}T$ and $\lambda p \in R(T^{*}T) \subset (N(T))^{\perp} \cap D(T^{*}T)$. It confirms that $\{p, 0\} \in (T^{*}T\vert_{(N(T))^{\perp}} - \lambda)$ (since $\lambda \neq 0$) and $\{0, p\} \in (T^{*}T\vert_{(N(T))^{\perp}} - \lambda)^{-1}$. The property of operator $(T^{*}T\vert_{(N(T))^{\perp}} - \lambda)^{-1}$ says that $p = 0$. So, $(T^{*}T - \lambda)^{-1}$ is an operator. $(T^{*}T - \lambda)^{-1}$ is closed because $(T^{*}T - \lambda)$ is closed. By Proposition II.6.3 \cite{MR1631548} and $R(T^{*}T\vert_{(N(T))^{\perp}} - \lambda) = (N(T))^{\perp}$, we get,
\begin{align*}
0 < \gamma(T^{*}T\vert_{(N(T))^{\perp}} - \lambda) = \gamma((T^{*}T - \lambda)\vert_{(N(T))^{\perp}}) \leq \gamma(T^{*}T - \lambda).
\end{align*}
From Proposition 3.1 \cite{MR3079830}, it is true that $D((T^{*}T - \lambda)^{-1}) = R(T^{*}T - \lambda)$ is closed. Thus, $(T^{*}T - \lambda)^{-1}$ is bounded. Again, the relation  $(T^{*}T\vert_{(N(T))^{\perp}} - \lambda)^{-1} \subset (T^{*}T - \lambda)^{-1}$ guarantees that $(N(T))^{\perp} \subset D((T^{*}T - \lambda)^{-1})$. Let $x_{0} \in N(T)$, then $\{x_{0}, 0\} \in T^{*}T$. So, $\{ -\lambda x_{0}, x_{0}\} \in (T^{*}T - \lambda)^{-1}$. Since, $\lambda \neq 0$ which implies $N(T) \subset D((T^{*}T - \lambda)^{-1})$. Thus, $D((T^{*}T - \lambda)^{-1}) = H$. Now, it is obvious to see that $\lambda \in \rho(T^{*}T)$ which is a contradiction. Our assumption is wrong. Moreover, 
\begin{align}\label{equ 3}
\sigma(T^{*}T) \setminus \{0\} \subset \sigma(T^{*}T\vert_{(N(T))^{\perp}}) \setminus \{0\}.
\end{align}
Now, the reverse inclusion will be shown. Consider $\mu \in \sigma(T^{*}T\vert_{(N(T))^{\perp}}) \setminus \{0\}$ but $\mu \in \rho(T^{*}T)$. Then, $(T^{*}T - \mu)^{-1}$ is a bounded operator in domain $H$. So, $(T^{*}T\vert_{(N(T))^{\perp}} - \mu)^{-1}$ is an bounded operator. Again, take an element $q \in (N(T))^{\perp} \subset H = R(T^{*}T - \mu)$. Then there exists $w \in D(T^{*}T)$ such that $\{w, q\} \in (T^{*}T - \mu)$. Again, $q + \mu w \in R(T^{*}T) \subset (N(T))^{\perp}$. Thus, $w \in (N(T))^{\perp}$ which implies $\{w, q\} \in (T^{*}T\vert_{(N(T))^{\perp}} - \mu)$ and $q \in D((T^{*}T\vert_{(N(T))^{\perp}} - \mu)^{-1})$. Hence, $(N(T))^{\perp} \subset  D((T^{*}T\vert_{(N(T))^{\perp}} - \mu)^{-1})$. Again, consider $s \in R(T^{*}T\vert_{(N(T))^{\perp}} - \mu)$, then there exists $t \in (N(T))^{\perp}$ such that $\{t, s\} \in (T^{*}T\vert_{(N(T))^{\perp}} - \mu)$. So, $s + \mu t \in (N(T))^{\perp}$ which implies $s \in (N(T))^{\perp}$. Furthermore, $(N(T))^{\perp} = D((T^{*}T\vert_{(N(T))^{\perp}} - \mu)^{-1})$. Now, it is ready to confirm that $\mu \in \rho(T^{*}T\vert_{(N(T))^{\perp}})$ which is again a contradiction. Therefore,
\begin{align}\label{equ 4}
\sigma(T^{*}T\vert_{(N(T))^{\perp}}) \setminus \{0\} \subset \sigma(T^{*}T) \setminus \{0\}.
\end{align}
By the relations (\ref{equ 3}) and (\ref{equ 4}), we have the relation
\begin{align}\label{equ 5}
\sigma(T^{*}T\vert_{(N(T))^{\perp}}) \setminus \{0\} = \sigma(T^{*}T) \setminus \{0\}.
\end{align}
\end{proof}
\begin{theorem}\label{thm 3.16}
Let $T \in CR(H, K)$. Then $\gamma(T^{*}T) = \gamma(T^{*}T\vert_{(N(T))^{\perp}})$.
\end{theorem}
\begin{proof}
By Lemma 5.1 \cite{MR2188974}, we get $T^{*}T$ is self-adjoint. Now, we claim that $T^{*}T\vert_{(N(T))^{\perp}}$ is self-adjoint. From Lemma \ref{lemma 3.14}, we have $T^{*}T\vert_{(N(T))^{\perp}}$ is closed. Moreover $T^{*}T\vert_{(N(T))^{\perp}}$ is symmetric because $T^{*}T\vert_{(N(T))^{\perp}} \subset T^{*}T \subset (T^{*}T\vert_{(N(T))^{\perp}})^{*}$. Theorem 1.5.5 \cite{MR3971207} says that $\lambda \in \rho(T^{*}T)$, when $\lambda \in \mathbb{C} \setminus \mathbb{R}$. By Theorem \ref{thm 3.15}, we get $\lambda \in \rho({T^{*}T\vert_{(N(T))^{\perp}}})$, for all $\lambda \in \mathbb{C} \setminus \mathbb{R}$. Again, Theorem 1.5.5 \cite{MR3971207}  confirms that $T^{*}T\vert_{(N(T))^{\perp}}$ is self-adjoint. From Theorem 4.3 \cite{MR3079830} and Theorem \ref{thm 3.15}, we can say that 
\begin{align*}
\gamma(T^{*}T) = \inf\{\vert \lambda \vert: \lambda \in \sigma(T^{*}T) \setminus \{0\}\} &=  \inf\{\vert \lambda \vert: \lambda \in \sigma(T^{*}T\vert_{(N(T))^{\perp}}) \setminus \{0\}\}\\ 
&= \gamma(T^{*}T\vert_{(N(T))^{\perp}}).
\end{align*}
\end{proof}
\begin{corollary}\label{cor 3.17}
Let $T \in CR(H, K)$. Then $T$ is Hyers-Ulam stable if and only if $T^{*}T\vert_{(N(T))^{\perp}}$ is Hyers-Ulam stable.
\end{corollary}
\begin{proof}
$T$ is Hyers-Ulam stable if and only if $R(T)$ is closed if and only if $R(T^{*})$ is closed if and only if $R(T^{*}T)$ is closed if and only if $\gamma(T^{*}T) = \gamma(T^{*}T\vert_{(N(T))^{\perp}}) > 0$ if and only if $R(T^{*}T\vert_{(N(T))^{\perp}})$ is closed if and only if $T^{*}T\vert_{(N(T))^{\perp}}$ is Hyers-Ulam stable.
\end{proof}
Now, we will discuss the sufficient conditions under which the sum and product of two Hyers-Ulam stable linear relations remain Hyers-Ulam stable.
\begin{theorem}\label{thm 3.18}
Let $T \in CR(H, K)$ and $S \in LR(H, K)$ such that $M(S) \subset M(T)$ and $D(T) \subset D(S)$ with the condition $\|Sx\| \leq b \|Tx\|$, for all $x \in D(T)$ and $0 < b < 1$. If $T$ is Hyers-Ulam stable, then $S+T$ is also Hyers-Ulam stable.
\end{theorem}
\begin{proof}
By Theorem 3.1.1 \cite{MR4492807}, we get $S + T$ is closed. It is enough to show that $R(S +T)$ is closed to get the Hyers-Ulam stability of $S + T$. It is obvious that $M(S + T) = M(T)$ because of the given condition $M(S) \subset M(T)$. Now, for all $x \in D(T)$, we have 
\begin{align*}
\|Sx + Tx\| = \|(Sx + Tx) + M(T)\| \leq \|Sx + M(T)\| + \|Tx + M(T)\| &\leq \|Sx\| + \|Tx\|\\
&\leq (1 + b) \|Tx\|.
\end{align*}
By Proposition 2.5.7 \cite{MR4492807} and Proposition II.1.5  \cite{MR1631548}, we get,
\begin{align*}
(1-b)\|Tx\| \leq \|Tx\|-\|Sx\| \leq \|(S + T)x\|, \text{ for all } x \in D(T).
\end{align*}
Combining the above two relations we get,
\begin{align}\label{equ 6}
(1-b)\|Tx\| \leq \|(S + T)x\| \leq (1 + b) \|Tx\|, \text{ for all } x \in D(T).
\end{align}
Now, we claim that $N(S + T) = N(T)$. Let $z \in N(S+ T)$, then $\{z, 0\} \in S +T$ such that $\{z, 0\} = \{z, z_{1}\} + \{z, z_{2}\}$, where $\{z, z_{1}\} \in S$ and $\{z, z_{2}\} \in T$. So, By the relation (\ref{equ 6}), we can say that
\begin{align*}
\|z_{2} + M(T) \| = \|Tz\| \leq \frac{1}{1-b} \|(S + T)z\| = 0.
\end{align*}
Thus, $z_{2} \in M(T)$ and $\{0, z_{2}\} \in T$ implies $\{z, 0\} \in T$. So, $N(S + T) \subset N(T)$. Again, consider $\{w, 0\} \in T$, then By the relation (\ref{equ 6}), we get $\|(S+T)w + M(S+T)\| = \|(S+T)w\| \leq (1 + b)\|Tw\| = 0.$ Thus, $(S+ T)w \subset M(S + T)$. Moreover, $\{w, 0\} = \{w, v_{1}\} + \{w, v_{2}\} \in (T +S) - S$ (by Proposition 2.3.4 \cite{MR4492807}), where $\{w, v_{1}\} \in T + S$ and $\{w, v_{2}\} \in -S$. So, $\{0, v_{1}\} \in S + T$ which implies $\{w, 0\} \in S + T$. Hence, $N(T) \subset N(S + T)$. Now, it is ready to say that $N(S + T) = N(T)$. Since $T$ is Hyers-Ulam stable. So, $T^{-1}$ is bounded because $T^{-1}$ is closed with $D(T^{-1}) = R(T)$ is closed. Again, $(S + T)^{-1}$ is closed because $S + T$ is closed.  We will show that $(S + T)^{-1}$ is bounded. Now, consider $\{q, p\} \in (S + T)^{-1}$, then $\{p, q\} =\{p, q^{'}\} + \{p, q^{''}\} \in S + T$, where $\{p, q^{'}\} \in S$ and $\{p, q^{''}\} \in T$. Moreover,
\begin{align}\label{equ 7}
\|(S + T)^{-1}q\| = \|(S + T)^{-1}q + M((S + T)^{-1})\| = \|(S + T)^{-1}q + N( T)\|.
\end{align}
The boundedness of $T^{-1}$  and the relation (\ref{equ 6}) say that
\begin{align}\label{equ 8}
(1-b) \|q^{''} + M(T)\| = (1-b) \|Tp\| \leq \|(S + T)p\| = \|q + M(S + T)\| \leq \|q^{'} + q^{''}\|.
\end{align}
Then there exists $s^{''} \in M(T)$ such that $\|q^{''} + s^{''}\| \leq \frac{1}{1-b}\|q^{'} + q^{''}\|.$ This confirms that $\{p, q^{''} + s^{''}\} \in T$. From  the relation(\ref{equ 7}) and (\ref{equ 8}), we have 
$$\|(S + T)^{-1}q\| = \|p + N(T)\| = \|T^{-1}(q^{''} + s^{''})+ M(T^{-1})\| \leq \|T^{-1}\|\|q^{''} + s^{''}\|  \leq \frac{\|T^{-1}\|}{1-b} \|q\|.$$
Hence, $(S + T)^{-1}$ is bounded. Furthermore, $R(S + T) = D((S + T)^{-1})$ is closed. Therefore, $S +T$ is Hyers-Ulam stable.
\end{proof}
Let us consider two linear relations $T \in LR(H_{1}, K_{1})$ and $S \in LR(H_{2}, K_{2})$, where $H_{i}$ and $K_{i}$ ($i = 1, 2$) are Hilbert spaces. We define the product of $T$ and $S$ by 
$$T \times S = \{\{(h_{1}, h_{2}), (k_{1}, k_{2})\}: \{h_{1}, k_{1}\} \in T \text{ and } \{h_{2}, k_{2}\} \in S\}.$$
It is easy to show that $T \times S$ is a linear relation from the Hilbert space $H_{1} \times H_{2}$ into the Hilbert space $K_{1} \times K_{2}$.
\begin{theorem}\label{thm 3.19}
Let $T \in CR(H_{1}, K_{1})$ and $S \in CR(H_{2}, K_{2})$ both be Hyers-Ulam stable. Then $T \times S$ is also Hyers-Ulam stable.
\end{theorem}
\begin{proof}
Let $\{x, y\} \in \overline{T \times S}$ (where $x \in H_{1} \times H_{2}$ and $y \in K_{1} \times k_{2}$), then there exists a sequence $\{\{(h_{1n}, h_{2n}), (k_{1n}, k_{2n})\}\}$ such that $\{\{(h_{1n}, h_{2n}), (k_{1n}, k_{2n})\}\} \to \{x, y\}$, as $n \to \infty$, where $\{h_{1n}, k_{1n}\} \in T$ and $\{h_{2n}, k_{2n} \} \in S$ for all $n \in \mathbb{N}$. Thus, $(h_{1n}, h_{2n}) \to x$ and $(k_{1n}, k_{2n}) \to y$, as $n \to \infty$. So, $\{h_{1n}\}$, $\{h_{2n}\}$, $\{k_{1n}\}$ and $\{k_{2n}\}$ are Cauchy sequences. We get some $h_{1} \in H_{1}$, $h_{2} \in H_{2}$, $k_{1} \in K_{1}$ and $k_{2} \in K_{2}$ such that $\{h_{1n}\} \to h_{1}$, $\{h_{2n}\} \to h_{2}$, $\{k_{1n}\} \to k_{1}$ and $\{k_{2n}\} \to k_{2}$ as $n \to \infty$. Moreover, $x = (h_{1}, h_{2})$ and $y= (k_{1}, k_{2})$ with $\{h_{1}, k_{1}\} \in T$ and $\{h_{2}, k_{2}\} \in S$ because $T$ and $S$ both are closed. Thus, $\{x, y\} \in T \times S$. Hence, $T \times S$ is closed. Again, consider $(y_{1}, y_{2}) \in \overline{R(T \times S)}$, there exists a sequence  $\{\{(h_{1n}^{'}, h_{2n}^{'}), (k_{1n}^{'}, k_{2n}^{'})\}\}$  in $T \times S$ such that $(k_{1n}^{'}, k_{2n}^{'}) \to (y_{1}, y_{2})$, as $n \to \infty$. $R(T)$ and $R(S)$ are closed because $T$ and $S$ both are Hyers-Ulam stable. Then $y_{1} \in R(T)$ and $y_{2} \in R(S)$. There exists $x_{1} \in D(T)$ and $x_{2} \in D(S)$ such that $\{x_{1}, y_{1}\} \in T$ and $\{x_{2}, y_{2}\} \in S$ which implies $(y_{1}, y_{2}) \in R(T \times S)$. Furthermore, $R(T \times S)$ is closed. Therefore, $T \times S$ is Hyers-Ulam stable.
\end{proof}
Theorem \ref{thm 3.20} depicts the Hyers-Ulam stability of a matrix linear relation. We define the matrix relation $\mathcal{A} = \begin{bmatrix}
  A & B\\
  C & F
  \end{bmatrix}$ from domain $(D(A) \cap D(C)) \times (D(B) \cap D(F)) \subset H \times K$ to $H \times K$ by
$$\mathcal{A} = \{\{(x,y),  (x_{a} + y_{b} , x_{c} + y_{f})\}: \{x, x_{a}\} \in A, \{x, x_{c}\} \in C, \{y, y_{b}\} \in B \text{ and } \{y, y_{f}\} \in F\}.$$, where $A \in LR(H, H)$, $B \in LR(K, H)$, $C \in LR(H, K)$ and $F \in LR(K, K)$ respectively. It is easy to show that $\mathcal{A}$ is a linear relation.
\begin{theorem}\label{thm 3.20}
Let  $\mathcal{A} = \begin{bmatrix}
  A & B\\
  C & F
  \end{bmatrix}$ be a matrix linear relation, where $A \in CR(H, H)$, $B \in LR(K, H)$, $C \in LR(H, K)$ and $F \in CR(K, K)$ respectively. Assume $M(B) \subset M(A)$ and $M(C) \subset M(F)$ with $\mathcal{A}$ is diagonally dominated which means $\|Cx\| \leq a \|Ax\|, \text{ for all } x \in D(A) \subset D(C)$ and $\|Bz\| \leq f \|Fz\|, \text{ for all } z \in D(F) \subset D(B)$, where $0 < a , f < 1$. Then $\mathcal{A}$ is a closed linear relation. Moreover, $\mathcal{A}$ is Hyers-Ulam stable when $A$ and $F$ both are Hyers-Ulam stable.
\end{theorem}
\begin{proof}
Let us define $T = \begin{bmatrix}
  A & 0\\
  0 & F
  \end{bmatrix}$ and $S = \begin{bmatrix}
  0 & B\\
  C & 0
  \end{bmatrix}$, where $\{(x, y), (x_{a}, y_{f})\} \in T$ and $\{(x, y), (y_{b}, x_{c})\} \in S$ for $\{x, x_{a}\} \in A$, $\{y, y_{f}\} \in F$, $\{y, y_{b}\} \in B$ and $\{x, x_{c}\} \in C$ with $x \in D(A)$ and $y \in D(F)$. By Theorem \ref{thm 3.19}, we get $T$ is a closed linear relation from $H \times K$ into $H \times K$. Again, $S$ is a linear relation from $H \times K$ into $H \times K$. Moreover, $D(T) \subset D(S)$. It is easy to show that
  $$M(S) = M(B) \times M(C) \subset M(A) \times M(F) = M(T).$$
  For all $(x, y) \in D(T)$, we get 
\begin{align*}
  \|S(x, y)\|^{2} &= \|(y_{b}, x_{c}) + (M (B) \times M(C))\|^{2}\\
 &= \|y_{b}+ M(B)\|^{2} + \|x_{c} + M(C)\|^{2} \\
 &= \|By\|^{2} + \|Cx\|^{2}\\
 & \leq a^{2}\|Ax\|^{2} + f^{2}\|Fy\|^{2}\\
 & \leq d^{2}\|T(x, y)\|^{2}, \text{ where }, 0< d = \max\{a, f\} < 1.
 \end{align*}
 Thus, $\|S(x, y)\| \leq d \|T(x, y)\|$, for all $(x, y) \in D(T)$. By Theorem 3.1.1 \cite{MR4492807}, we can say that $\mathcal{A} = T + S$ is closed. By Theorem \ref{thm 3.19}, we have $T$ is Hyers-Ulam stable because $A$ and $F$ both are Hyers-Ulam stable. Therefore, Theorem \ref{thm 3.18} confirms that $\mathcal{A} = T + S$ is also Hyers-Ulam stable.
\end{proof}

\section{Conclusions}
 In this paper, the Hyers-Ulam stability of linear relations in normed linear spaces is introduced, and several interesting results concerning the Hyers-Ulam stability of closed linear relations in Hilbert spaces are explored. By establishing the result $\sigma(T^{*}T\vert_{(N(T))^{\perp}}) \setminus \{0\} = \sigma(T^{*}T) \setminus \{0\}$ (when $T$ is a closed relation from Hilbert space $H$ into $K$), it is proved that $T$ is Hyers-Ulam stable if and only if $T^{*}T\vert_{(N(T))^{\perp}}$ is Hyers-Ulam stable. Additionally, sufficient conditions are provided under which the sum and product of two Hyers-Ulam stable linear relations remain Hyers-Ulam stable.

\section*{Declarations}
  The author declares that there are no conflicts of interest.


\begin{thebibliography}{10}
   
\bibitem{majumdar2024hyers}
 Arup Majumdar, P Sam Johnson, and Ram N Mohapatra.\newblock{\em Hyers-Ulam Stability of Unbounded Closable Operators in Hilbert Spaces}. \newblock{Mathematische Nachrichten}, 297:10, 3887-3903, 2024.


  \bibitem{MR007}
   Arup Majumdar, P. Sam Johnson.
   \newblock{The Moore-Penrose inverses of unbounded closable operators and the direct sum of closed operators in Hilbert spaces}, \newblock{\em Linear and Multilinear Algebra}, 2024.

\bibitem{MRARUP}
Arup Majumdar, P. Sam Johnson.
\newblock{Characterizations of closed EP operators on Hilbert spaces}.
\newblock{\em arXiv:2410.06869v2}, 2024.
    
\bibitem{MR4492807}
Aymen Ammar and Aref Jeribi.
 \newblock{\em Spectral theory of multivalued linear operators},
 \newblock{Apple Academic Press, Oakville, ON; CRC Press, Boca Raton, FL}, 2022.


\bibitem{MR1427262}
   Balmohan V. Limaye.
     \newblock{\em Functional analysis}, \newblock{New Age International Publishers Limited, New Delhi}, 612, 1996.
      


\bibitem{MR1671909}
Claudi Alsina and Roman Ger.
\newblock {\emph {On some inequalities and stability results related to the exponential
function}}.
\newblock {J. Inequal. Appl.}, 2(4):373--380, 1998.
 


 \bibitem{MR2204863}
	Go~Hirasawa and Takeshi Miura.
	\newblock {\emph{Hyers-{U}lam stability of a closed operator in a {H}ilbert space}}.
	\newblock{Bull. Korean Math. Soc.}, 43(1):107--117, 2006.
 
\bibitem{MR888141}
   Gyokai Nikaido.
    \newblock{\em Product of linear operators with closed range},
   \newblock{Proc. Japan Acad. Ser. A Math. Sci.}, 62(9):338--340, 1986.
 

\bibitem{MR3971207}
Jussi Behrndt, Seppo Hassi, and Henk de Snoo.
     \newblock{\em Boundary value problems, {W}eyl functions, and differential
              operators}. \newblock{{Birkh\"auser/Springer, Cham}}, (108):772, 2020.
    




    \bibitem{MR1321558}
	Marta Ob{\l}oza.
	\newblock{\emph{Hyers stability of the linear differential equation}}.
	\newblock {Rocznik Nauk.-Dydakt. Prace Mat.}, (13):259--270, 1993.
   

\bibitem{MR1631548}
    Ronald Cross.
     \newblock{\em Multivalued linear operators}.
     \newblock{Marcel Dekker, Inc., New York}, (213):335, 1998.
    
\bibitem{MR2188974}
   S. Hassi, A. Sandovici, H. S. V. De Snoo,  and H. Winkler.
   \newblock{\em Form sums of nonnegative selfadjoint operators}, \newblock{Acta Math. Hungar.}, 111:(1-2):81--105, 2006.
  
    
   
   
\bibitem{MR2000046}
	Takeshi Miura, Shizuo Miyajima, and Sin-Ei Takahasi.
	\newblock {\emph{Hyers-{U}lam stability of linear differential operator with constant coefficients}}.
	\newblock {Math. Nachr.}, 258:90--96, 2003.
	


\bibitem{MR3079830}
   Teresa Alvarez and Adrian Sandovici.
     \newblock{\em On the reduced minimum modulus of a linear relation in
              {H}ilbert spaces}, \newblock{Complex Anal. Oper. Theory}, 7(4):801--812, 2013.

\bibitem{MR1778016}
	Themistocles~M. Rassias.
	\newblock {\emph{On the stability of functional equations and a problem of {U}lam}}.
	\newblock {Acta Appl. Math.}, 62(1):23--130, 2000.
	
 
  \bibitem{MR4203636}
  Zsigmond Tarcsay and Zolt\'an Sebesty\'en.
     \newblock{\em Canonical graph contractions of linear relations on {H}ilbert
              spaces}, \newblock{Complex Anal. Oper. Theory}, 15(1):16, 2021.
  
  
  

\end{thebibliography}
\end{document}